\documentclass[12pt]{article}
\usepackage{latexsym,amsfonts,amssymb}
\usepackage[dvips]{color}
\usepackage{colordvi,multicol}

\newtheorem{theorem}{Theorem}[section]
\newtheorem{corollary}[theorem]{Corollary}
\newtheorem{lemma}[theorem]{Lemma}
\newtheorem{proposition}[theorem]{Proposition}
\newtheorem{definition}[theorem]{Definition}
\newtheorem{remark}[theorem]{Remark}

\begin{document}

\title{\bf Multi-dimensional central limit theorems and  laws of
large numbers under sublinear expectations}
\author{Ze-Chun Hu\thanks{Corresponding
author: Department of Mathematics, Nanjing University, Nanjing
210093, China\vskip 0cm E-mail address: huzc@nju.edu.cn}\quad and Ling Zhou \\
 {\small Nanjing University}}
 \date{}
 \maketitle

\noindent{\bf Abstract}\quad In this paper, we present some multi-dimensional central limit theorems and  laws of large numbers under sublinear expectations, which extend some
previous results.


\smallskip

\noindent {\bf Keywords}\quad central limit theorem, laws of large numbers, sublinear expectation

\noindent {\bf MSC(2000):}  60H10, 60H05

\smallskip



\section{Introduction}
The classic central limit theorems and strong (weak) laws of large numbers  play an important role in the development of probability theory and its applications.
Recently, motivated by the risk measures,  superhedge pricing and modeling uncertain in finance, Peng \cite{{Peng1,Peng2,Peng3,Peng4,Peng5,Peng6}}  initiated the notion of independent and identically distributed (IID) random  variables under sublinear expectations, and proved the central limit theorems  and the weak law of large  numbers among other things.

 In \cite{HuZhang}, Hu and Zhang obtained a central limit theorem for capacities induced by sublinear expectations.  In \cite{LS}, Li and Shi proved a  central limit theorem without the requirement of identical distribution, which extends Peng's central limit theorems in \cite{Peng2,Peng4}. In \cite{Hu11}, Hu proved that for any continuous function $\psi$ satisfying the growth condition $|\psi(x)|\leq C(1+|x|^p)$ for some $C>0,p\geq 1$ depending on $\psi$, central limit theorem under sublinear expectations obtained by Peng \cite{Peng4} still holds. In \cite{Zhang}, Zhang obtained a central limit theorem for weighted sums of independent random variables under sublinear expectations, which extends the results obtained by Peng, Li and Shi.

 Now we recall the existing results about the laws of large numbers under sublinear expectations.  In \cite{Chen10}, Chen proved a
 strong law of large numbers for IID random variables
 under capacities induced by sublinear  expectations. In \cite{HuF}, Hu
 presented three laws of large numbers for independent random
 variables without the requirement of identical distribution. In \cite{HuYang}, we obtained  a strong law of large numbers for IID random variables under one-order type moment condition.

 The main results in \cite{HuZhang,LS,Hu11,Zhang,Chen10,HuF, HuYang} are about one-dimensional random variables in sublinear expectation space.
 The goal of this paper is to present some multi-dimensional central limit theorems and laws of large numbers under sublinear expectations without the requirement of identical distribution, which extend  Peng's central limit theorems in \cite{Peng2,Peng4}, Theorem 3.1 of \cite{LS} and Theorem 3.1 of \cite{HuF}.

The rest of this paper is organized as follows. In Section 2, we recall some basic notions under sublinear
expectations. In Section 3, we present a multi-dimensional central limit theorem and several corollaries, and their proofs will be given  in Section 4. In Section 5, we give some multi-dimensional  laws of large numbers.

\section{Basic settings}

In this section, we present some basic settings about sublinear expectations. Please refer to  Peng \cite{Peng1,Peng2,Peng3,Peng4,Peng5,Peng6}.

Let $\Omega$ be a given set and let $\bf{\mathcal {H}}$ be a linear space of real functions defined on $\Omega$ such that
 for any constant number $c, c\in \mathcal{H}$;  if $X\in\mathcal{H}$, then $|X|\in\mathcal{H}$; if $X_1,\ldots,X_n\in\mathcal{H}$, then for any $\varphi\in C_{l,Lip}(\mathbb{R}^n)$, $\varphi(X_1,\ldots,X_n)\in\mathcal{H}$, where $C_{l,Lip}(\mathbb{R}^n)$ denotes the linear space of functions $\varphi$ satisfying
$$|\varphi(x)-\varphi(y)|\le C(1+|x|^m+|y|^m)|x-y|, \forall x,y\in \mathbb{R}^n,$$
for some $C>0, m\in \mathbb{N}$ depending on $\varphi$.
Denote by $C_{b,Lip}(\mathbb{R}^n)$  the space of bounded Lipschitz
functions defined on $\mathbb{R}^n$.

\begin{definition}
A  sublinear expectation ${\hat{E}}$ on $H$ is a functional $\hat{E}:\mathcal{H}\rightarrow\mathbb{R}$ satisfying the
following properties:\\
\hspace*{0.3cm} (a) Monotonicity: $\hat{E}[X]\ge\hat{E}[Y]$, if $X\ge Y$.\\
\hspace*{0.3cm}  (b) Constant preserving: $\hat{E}[c]=c,\forall c\in\mathbb{R}.$\\
\hspace*{0.3cm}  (c) Sub-additivity: $\hat{E}[X+Y]\le\hat{E}[X]+\hat{E}[Y].$\\
\hspace*{0.3cm} (d) Positive homogeneity: $\hat{E}[\lambda X]=\lambda\hat{E}[X]$, $\forall\lambda\ge 0.$\\
The triple $(\Omega,\mathcal{H},\hat{E})$ is called a
sublinear expectation space. If only (c) and (d) are satisfied, $\hat{E}$ is called a sublinear functional.
\end{definition}

\begin{theorem}\label{Th1}
Let $\mathbb{E}$ be a sublinear functional defined on a linear space $\mathcal{H}$. Then there exists a family of linear functionals
$\{\mathrm{E}_{\theta}: \theta\in\Theta\}$ defined on $\mathcal{H}$ such that
$$
\mathbb{{E}}[X]=\sup_{\theta\in\Theta}\mathrm{E}_{\theta}[X],
\ \forall X\in\mathcal{H},
$$
and, for each $X\in\mathcal{H}$, there exists $\theta_X\in\Theta$ such that $\mathbb{E}[X]=\mathrm{E}_{\theta_X}[X]$. Furthermore, if $\mathbb{{E}}$ is a sublinear expectation, then the corresponding $\mathrm{E}_{\theta}$ is a linear expectation.
\end{theorem}

\begin{definition}
Let $X_1$ and $X_2$ be two $n$-dimensional random vectors defined on
sublinear expectation spaces $(\Omega_1,\mathcal{H}_1,E_1)$ and $(\Omega_2,\mathcal{H}_2,E_2)$, respectively. They
are called identically distributed, denoted by ${X_1\thicksim X_2}$, if
$$
E_1[\varphi(X_1)]=E_2[\varphi(X_2)],\ \forall\varphi\in C_{b,Lip}({\mathbb{R}}^n).
$$
\end{definition}

\begin{definition}
In a sublinear expectation space $(\Omega,\mathcal{H},\hat{E})$, a random vector
$Y\in{\mathcal{H}}^n$ is said to be independent from another random vector $X\in{\mathcal{H}}^m$ under $\hat{E}[\cdot]$ if for each test function $\varphi \in C_{l,Lip}{(\mathbb{R}^{n+m})}$ we have
$$\hat{E}[\varphi(X,Y)]=\hat{E}[\hat{E}[\varphi(x,Y)]|_{x=X}].$$
\end{definition}


\begin{proposition}
Let $(\Omega,\mathcal{H},\hat{E})$ be a sublinear expectation space and $X,Y$ be
two random variables such that $\hat{E}[Y]=-\hat{E}[-Y]$, i.e., $Y$ has no mean-uncertainty. Then we have
$$\hat{E}[X+\alpha Y]=\hat{E}[X]+\alpha \hat{E}[Y], \forall\alpha \in \mathbb{R}.$$
In particular, if $\hat{E}[Y]=\hat{E}[-Y]=0,$ then $\hat{E}[X+\alpha Y]=\hat{E}[X].$
\end{proposition}

\begin{proposition}
For each $X,Y\in\mathcal{H}$, we have
\begin{eqnarray*}
&&\hat{E}[|X+Y|^r]\le 2^{r-1}(\hat{E}[|X|^r]+\hat{E}[|Y|^r]),\\
&&\hat{E}[|XY|]\le (\hat{E}[|X|^p])^{\frac{1}{p}}\cdot(\hat{E}[|X|^q])^{\frac{1}{q}},\\
&&(\hat{E}[|X+Y|^p])^{\frac{1}{p}}\le (\hat{E}[|X|^p])^{\frac{1}{p}}+(\hat{E}[|Y|^p])^{\frac{1}{p}},
\end{eqnarray*}
where $r\ge 1$ and $1<p,q<\infty$ with $\frac{1}{p}+\frac{1}{q}=1.$
In particular, for $1\le p<q$, we have $$(\hat{E}[|X|^p])^{\frac{1}{p}}\le (\hat{E}[|X|^q])^{\frac{1}{q}}.$$
\end{proposition}

\begin{definition}
Let $X,\bar{X}$ be two n-dimensional random vectors on a sublinear expectation space $(\Omega,\mathcal{H},\hat{E})$. $\bar{X}$ is called an independent copy of $X$ if $\bar{X}\thicksim X$ and $\bar{X}$ is independent from $X$.
\end{definition}

\begin{definition}
($G$-normal distribution) A $d$-dimension random
variable $X=(X_1,X_2,\ldots,X_d)^T$ on a sublinear expectation space $(\Omega,\mathcal{H},\hat{E})$ is called (centralized) $G$-normal distributed if
$$aX+b\bar{X}\thicksim\sqrt{a^2+b^2}X,\ \mbox{for}\ a,b\ge 0,$$
where $\bar{X}$ is an independent copy of $X$.
\end{definition}

We denote by $\mathbb{S}(d)$ the collection of all $d\times d$ symmetric matrices. Let $X$
be $G$-normal distributed on $(\Omega,\mathcal{H},\hat{E})$. The following function is very important to characterize
its distribution:
$$
G(A)=\hat{E}\left[\frac{1}{2}\langle AX,X\rangle\right],~A\in\mathbb{S}(d).
$$
It is easy to check that $G$ is a sublinear function monotonic in $A\in \mathbb{S}(d)$. By Theorem \ref{Th1}, there
exists a bounded, convex and closed subset $\Theta\subset \mathbb{S}(d)$ such that
$$G(A)=\sup_{Q\in \Theta}\frac{1}{2}{\rm tr}[AQ],~A\in \mathbb{S}(d).$$

\begin{proposition}\label{Pro1}
Let $X$ be $G$-normal distributed. For each $\varphi\in C_{b,Lip}(\mathbb{R}^d)$, define a function $$u(t,x):=\hat{E}[\varphi(x+\sqrt{t}X)], ~\forall (t,x)\in [0,\infty)\times \mathbb{R}^d.$$
Then $u$ is the unique viscosity solution of the following parabolic PDE:
\begin{equation}\label{eq18}
\partial _t u-\textit{G}(D^2u)=0,~u|_{t=0}=\varphi,
\end{equation}
where $G(A)=\hat{E}[\frac{1}{2}\langle AX,X\rangle],~ A\in\mathbb{S}(d).$
\end{proposition}

The parabolic PDE (\ref{eq18}) is called a $G$-heat equation.
Since $G(A)$ is monotonic: $G(A_1) \ge G(A_2)$, for $A_1 \ge A_2$, it follows that
$$\Theta\subset \mathbb{S}_{+}(d)=\{\theta\in\mathbb{S}(d):\theta\ge 0\}=\{BB^T:B\in\mathbb{R}^{d\times d}\}.$$
If $\Theta=\{Q\}$, then $X$ is classical zero-mean normal distributed with covariance $Q$. In general, $\Theta$
characterizes the covariance uncertainty of $X$. We denote $X\thicksim N(0;\Theta)$.

When $d=1$, we have $X\thicksim N(0;[\underline{\sigma}^2,\overline{\sigma}^2])$, where $\overline{\sigma}^2=\hat{E}[X^2],~\underline{\sigma}^2=-\hat{E}[-X^2]$. The corresponding $G$-heat equation is
\begin{equation}\label{eq19}
\partial_t u-\frac{1}{2}(\overline{\sigma}^2(\partial_{xx}^2 u)^{+}-\underline{\sigma}^2(\partial_{xx}^2 u)^{-})=0,~~u|_{t=0}=\varphi.
\end{equation}

\section{Multi-dimensional central limit theorems under sublinear expectations}\setcounter{equation}{0}

 In this section, we present a multi-dimensional central limit theorem under sublinear expectation and several corollaries. Their proofs will be given in next section.

\begin{theorem}\label{Thm3.1}
Let  $\{X_i\}_{i=1}^{\infty}$ be a sequence of $\ \mathbb{R}^d$-valued random variables
in a sublinear expectation space $(\Omega,\mathcal{H},\hat{E})$ which satisfies the following conditions:
\begin{itemize}
\item[(i)] each $X_{i+1}$ is independent from $(X_1,\ldots,X_i),\forall i=1,2,\ldots;$
\item[(ii)] $\hat{E}[X_i]=\hat{E}[-X_i]=0;$
\item[(iii)] there is a constant $M>0$, such that $\hat{E}[|X_i|^3]\le M, \forall i=1,2,\ldots;$
\item[(iv)] there exist a  sequence $\{a_n\}$ of positive  numbers and a sublinear function $G:\mathbb{S}(d)\rightarrow\mathbb{R}$ such that
$\lim_{n\to \infty}\frac{a_1+\cdots+a_n}{n}=0$ and
\begin{eqnarray*}
|G_n(A)- G(A)|\le a_n\|A\|,~\forall n, A\in \mathbb{S}(d),
\end{eqnarray*}
where $G_i(A)=\frac{1}{2}\hat{E}[\langle AX_i,X_i\rangle],\|A\|=\sqrt{\sum_{i,j=1}^na^2_{ij}},\ \forall A=(a_{ij})\in \mathbb{S}(d)$.
\end{itemize}
Then the sequence $\left\{\frac{S_n}{\sqrt{n}}\right\}_{n=1}^{\infty}$ converges in law to $G$-normal distribution $N(0;\Theta)$, i.e.,
\begin{equation}\label{eq10}
\lim_{n\rightarrow\infty}\hat{E}\left[\varphi\left(\frac{S_n}{\sqrt{n}}
\right)\right]=\hat{E}[\varphi(X)],\ \forall \varphi\in C_{b,Lip}(\mathbb{R}^d),
\end{equation}
where $X\thicksim N(0;\Theta), S_n=X_1+\cdots+X_n$.
\end{theorem}

\bigskip

\begin{corollary}\label{Cor3.2} Let  $\{X_i\}_{i=1}^{\infty}$ be a sequence of $\ \mathbb{R}^d$-valued random variables
in a sublinear expectation space $(\Omega,\mathcal{H},\hat{E})$ which satisfies conditions (i)(ii)(iii) in  Theorem \ref{Thm3.1} and  the following condition (v).
\begin{itemize}
\item[(v)] There exists a sublinear function $G: \mathbb{S}(d)\rightarrow\mathbb{R}$ such that $\{G_i\}$  converges to $G$ pointwise, i.e.,
$$\lim_{i \to \infty } G_i(A)=G(A),\ \forall A\in \mathbb{S}(d),$$
where $G_i(A)=\frac{1}{2}\hat{E}[\langle AX_i,X_i\rangle],\forall A\in \mathbb{S}(d)$.
\end{itemize}
Then the result of Theorem \ref{Thm3.1} holds.
\end{corollary}
\bigskip

In order to state another corollary of Theorem \ref{Thm3.1},  we need the following definition.

\begin{definition}\label{def3.3}
Let $\Delta _1,\Delta_2\subset \mathbb{S}_+(d)$ be two closed subsets. Define their Hausdorff distance as follows:
$$
d_H(\Delta_1,\Delta_2):=\inf\{\varepsilon |\Delta_1\subset B_\varepsilon(\Delta_2),~\Delta_2\subset B_\varepsilon(\Delta_1)\},$$
where
$B_\varepsilon(\Delta):= \{x\in\mathbb{S}_+(d)|d(x,\Delta)<\varepsilon\},
\ d(x,\Delta):= \inf\{d(x,y)|y\in \Delta\},$ and
$$d(x,y):= \|x-y\|=\sqrt{\sum_{i,j=1}^d|x_{ij}-y_{ij}|^2}\ \mbox{for}\ x=(x_{ij}),y=(y_{ij}).$$
\end{definition}

Set $\Xi:= \{\Delta\subset\mathbb{S}_+(d)|\Delta\ \mbox{is bounded and closed}\}.$
Then by \cite{He},  $(\Xi,d_H)$ is a complete metric space.

\bigskip

\begin{corollary}\label{cor3.4} Let  $\{X_i\}_{i=1}^{\infty}$ be a sequence of $\ \mathbb{R}^d$-valued random variables
in a sublinear expectation space $(\Omega,\mathcal{H},\hat{E})$ which satisfies conditions (i)(ii)(iii) in  Theorem \ref{Thm3.1} and  the following condition (vi).
\begin{itemize}
\item[(vi)] There exists  a bounded, convex, and  closed subset $\Theta$ of  $\ \mathbb{S}_+(d)$ such that $$\lim_{n\rightarrow \infty}\frac{1}{n}\sum_{i=1}^{n}d_H(\Theta,\Theta_i)=0,$$
    where  $\Theta_i$ is a bounded, convex and   closed  subset of $\ \mathbb{S}_+(d)$  such that
$$G_i(A):=\frac{1}{2}\hat{E}[\langle AX_i,X_i\rangle]=\frac{1}{2}\sup_{Q\in\Theta_i}\rm{tr}[AQ].$$
\end{itemize}
Define $G(A)=\frac{1}{2}\sup_{Q\in\Theta}{\rm tr}[AQ],\forall A\in \mathbb{S}(d).$
Then the result of Theorem \ref{Thm3.1} holds.
\end{corollary}

\bigskip

\begin{corollary}\label{Cor3.5} Let  $\{X_i\}_{i=1}^{\infty}$ be a sequence of $\ \mathbb{R}^d$-valued random variables
in a sublinear expectation space $(\Omega,\mathcal{H},\hat{E})$ which satisfies conditions (i)(ii)(iii) in  Theorem \ref{Thm3.1} and  the following condition (vii).
\begin{itemize}
\item[(vii)] Let $\Theta_i$ be a bounded, convex and closed subset of $\ \mathbb{S}_+(d)$  such that
$$G_i(A):=\frac{1}{2}\hat{E}[\langle AX_i,X_i\rangle]=\frac{1}{2}\sup_{Q\in\Theta_i}\rm{tr}[AQ],$$
and $\{\Theta _i\}$ is a Cauchy sequence in $(\Xi,d_H)$, i.e.,
$d_H(\Theta _i,\Theta _j)\to 0,\ \mbox{as}\ i,j \to \infty.$
\end{itemize}
Then there exists a bounded, convex and   closed   subset $\Theta\subset \mathbb{S}_+(d)$ such that the result of Theorem \ref{Thm3.1} holds, where the function $G: \mathbb{S}(d)\rightarrow \mathbb{R}$ is defined by   $G(A)=\frac{1}{2}\sup_{Q\in\Theta}{\rm tr}[AQ].$
\end{corollary}

\bigskip

By Corollary \ref{cor3.4}, we can obtain the following result.

\begin{theorem}(\cite{LS})\label{thm3.6}
Let a sequence $\{X_i\}_{i=0}^{\infty}$ which is 1-dimensional random variable in a sublinear expectation space $(\Omega,\mathcal{H},\hat{E})$ satisfy the
following conditions:
\begin{itemize}
\item[(i)] each $X_{i+1}$ is independent from $(X_1,\ldots,X_i)$, $\forall i=1,2,\ldots;$
\item[(ii)] $\hat{E}[X_i]=\hat{E}[-X_i]=0$, $\hat{E}[X_i^2]=\overline{\sigma}_i^2$, $-\hat{E}[-X_i^2]=\underline{\sigma}_i^2$, where
$0\le\underline{\sigma}_i^2 \le \overline{\sigma}_i^2<\infty;$
\item[(iii)] there are two positive constants $\underline{\sigma}$ and $\overline{\sigma}$ such that
$$\lim_{n\rightarrow \infty}\frac{1}{n}\sum_{i=1}^{n}|\underline{\sigma}_i^2-\underline{\sigma}^2|=0,~~
\lim_{n\rightarrow \infty}\frac{1}{n}\sum_{i=1}^{n}|\overline{\sigma}_i^2-\overline{\sigma}^2|
=0;$$
\item[(iv)] there is a constant $M>0$, such that $$\hat{E}[|X_i|^3]\le M,~\forall i=1,2,\ldots$$
\end{itemize}
Then the sequence $\left\{\frac{S_n}{\sqrt{n}}\right\}_{n=1}^{\infty}$ of the sum $S_n=X_1+\cdots+X_n$ converges in law to  $N(0;[\underline{\sigma}^2,\overline{\sigma}^2]):$
$$\lim_{n\rightarrow\infty}\hat{E}\left[\varphi\left(\frac{S_n}
{\sqrt{n}}\right)\right]=\hat{E}[\varphi(X)],\ \forall \varphi\in C_{b,Lip}(\mathbb{R}),$$
where $X\thicksim N(0;[\underline{\sigma}^2,\overline{\sigma}^2])$.
\end{theorem}
{\bf Proof.}
For $d=1$, by Definition \ref{def3.3}, we have
$$
d_H([a,b],[c,d])=\max\{|a-c|,|b-d|\}\le |a-c|+|b-d|.
$$
Set $\Theta_i=[\underline{\sigma}_i^2 , \overline{\sigma}_i^2]$, $\Theta=[\underline{\sigma}^2 , \overline{\sigma}^2]$.
Then
$$
d_H(\Theta_i,\Theta)=\max\{|\underline{\sigma}_i^2-\underline{\sigma}^2| ,|\overline{\sigma}_i^2-\overline{\sigma}^2|\}
\le |\underline{\sigma}_i^2-\underline{\sigma}^2| +|\overline{\sigma}_i^2-\overline{\sigma}^2|,
$$
and thus by condition (iii), we have
$$
\lim_{n\rightarrow \infty}\frac{1}{n}\sum_{i=1}^{n}d_H(\Theta_i,\Theta)\le \lim_{n\rightarrow \infty}\frac{1}{n}\sum_{i=1}^{n}|\underline{\sigma}_i^2-\underline{\sigma}^2|+\lim_{n\rightarrow \infty}\frac{1}{n}\sum_{i=1}^{n}|\overline{\sigma}_i^2-
\overline{\sigma}^2|=0.
$$
Then the result follows from Corollary \ref{cor3.4}.
\hfill\fbox

\section{Proofs}\setcounter{equation}{0}
In this section, we will give the proofs of Theorem \ref{Thm3.1} and its corollaries.

\subsection{Proof of Theorem \ref{Thm3.1}}


At first, we prove a lemma.

\begin{lemma}\label{lem4.1}
We assume the same assumptions as in Theorem \ref{Thm3.1}. We further
assume that there exists a constant $\beta>0$ such that, for each $A,B\in \mathbb{S}(d)$ with $A \ge B$, we have
\begin{equation}\label{lem4.1-a}
\hat{E}[\langle AX,X\rangle]-\hat{E}[\langle BX,X\rangle]\ge \beta\ {\rm tr}[A-B].
\end{equation}
Then the main result $($\ref{eq10}$)$ holds.
\end{lemma}
{\bf Proof.}
 For a small but fixed $h>0$, let $V(t,x)$ be the unique viscosity solution of
\begin{equation}\label{lem4.1-b}
\partial_tV+G(D^2V)=0,~(t,x)\in[0,1+h]\times\mathbb{R}^d,~V|_{t=1+h}=\varphi.
\end{equation}
By Proposition \ref{Pro1}, we have
\begin{eqnarray}\label{lem4.1-b-1}
V(t,x)=\hat{E}[\varphi(x+\sqrt{1+h-t}X)].
\end{eqnarray}
In particular,
$$
V(h,0)=\hat{E}[\varphi(X)],~V(1+h,x)=\varphi(x).
$$
Since (\ref{lem4.1-b}) is a uniformly parabolic PDE and $G$ is a convex function, by the
interior regularity of $V$ (see \cite{Wang}), we have
\begin{eqnarray}\label{lem4.1-c}
\|V\|_{C^{1+\frac{\alpha}{2},2+\alpha}([0,1]\times\mathbb{R}^d)}<\infty, \ \mbox{for some}\ \alpha\in(0,1).
\end{eqnarray}

Set $\delta=\frac{1}{n}$, $S_0=0$. Then
\begin{eqnarray*}\label{eq5}
V(1,\sqrt{\delta}S_n)-V(0,0)
&=&\sum_{i=0}^{n-1}\left\{V((i+1)\delta,\sqrt{\delta}S_{i+1})-V(i\delta,\sqrt{\delta}S_i)\right\}\\
&=&\sum_{i=0}^{n-1}\{[V((i+1)\delta,\sqrt{\delta}S_{i+1})-V(i\delta,\sqrt{\delta}S_{i+1})]\\
&&+[V(i\delta,\sqrt{\delta}S_{i+1})-V(i\delta,\sqrt{\delta}S_i)]\}\\
&=&\sum_{i=0}^{n-1}\{I_\delta^i+J_\delta^i\},
\end{eqnarray*}
with, by Taylor¡¯s expansion,
\begin{eqnarray*}
J_\delta^i&=&\partial_tV(i\delta,\sqrt{\delta}S_i)\delta+
\frac{1}{2}\langle D^2V(i\delta,\sqrt{\delta}S_i)X_{i+1},X_{i+1}\rangle\delta
+\langle DV(i\delta,\sqrt{\delta}S_i),X_{i+1}\sqrt{\delta}\rangle,\\
I_\delta^i&=&\delta\int_0^1[\partial_tV((i+\beta)\delta,\sqrt{\delta}
S_{i+1})-\partial_tV(i\delta,\sqrt{\delta}S_{i+1})]\rm{d}\beta\\
           &&+[\partial_tV(i\delta,\sqrt{\delta}S_{i+1})-
           \partial_tV(i\delta,\sqrt{\delta}S_{i})]\delta\\
           &&+\int_0^1\int_0^1\delta\left\langle \left(D^2V(i\delta,\sqrt{\delta}S_{i}+\gamma\beta X_{i+1}\sqrt{\delta})-D^2V(i\delta,\sqrt{\delta}S_i)\right)X_{i+1} ,X_{i+1}\right\rangle \gamma \rm{d}{\beta}\rm{d}\gamma.
\end{eqnarray*}
Thus
\begin{eqnarray}\label{eq6}
\hat{E}\left[\sum_{i=0}^{n-1}J_\delta^i\right]-
\hat{E}\left[-\sum_{i=0}^{n-1}I_\delta^i\right]
&\le&\hat{E}[V(1,\sqrt{\delta}S_n)]-V(0,0)
\le\hat{E}\left[\sum_{i=0}^{n-1}J_\delta^i\right]+\hat{E}
\left[\sum_{i=0}^{n-1}I_\delta^i\right].\ \
\end{eqnarray}

We now prove that
$\lim_{n\rightarrow\infty}\hat{E}\left[\sum_{i=0}^{n-1}J_\delta^i\right]=0.
$
Firstly, by conditions (i) and (ii), we have
$$\hat{E}[\langle DV(i\delta,\sqrt{\delta}S_i),X_{i+1}\rangle \sqrt{\delta}]=\hat{E}[\langle -DV(i\delta,\sqrt{\delta}S_i),X_{i+1}\rangle \sqrt{\delta}]=0.$$
Then by (\ref{lem4.1-b}) and the sublinear property of $\hat{E}$, we have
\begin{eqnarray*}\label{eq7}
\hat{E}\left[\sum_{i=0}^{n-1}J_\delta^i\right]
&=&\hat{E}\left[\sum_{i=0}^{n-1}\left(\partial_tV(i\delta,\sqrt{\delta}S_i)
   \delta+\frac{1}{2}\langle D^2V(i\delta,\sqrt{\delta}S_i)X_{i+1},X_{i+1}\rangle\delta\right)\right]\\
&=&\hat{E}\left[\sum_{i=0}^{n-1}\left(\partial_tV(i\delta,\sqrt{\delta}S_i)
   \delta+G(D^2V(i\delta,\sqrt{\delta}S_i))\delta\right)\right.\\
&&+\left.\sum_{i=0}^{n-1}\left(\frac{1}{2}\langle
   D^2V(i\delta,\sqrt{\delta}S_i)X_{i+1},X_{i+1}\rangle\delta-
   G(D^2V(i\delta,\sqrt{\delta}S_i))\delta\right)\right]\\
&=&\hat{E}\left[\sum_{i=0}^{n-1}\left(\frac{1}{2}\langle
   D^2V(i\delta,\sqrt{\delta}S_i)X_{i+1},X_{i+1}\rangle-
   G(D^2V(i\delta,\sqrt{\delta}S_i))\right)\delta\right]\\
&\le&  \delta \hat{E}\left[\sum_{i=0}^{n-1}\left(\frac{1}{2}\langle D^2V(i\delta,\sqrt{\delta}S_i)X_{i+1},X_{i+1}\rangle-
G_{i+1}(D^2V(i\delta,\sqrt{\delta}S_i))\right)\right]\\
&&+\delta \hat{E}\left[\sum_{i=0}^{n-1}\left(G_{i+1}(D^2V(i\delta,\sqrt{\delta}S_i))-
G(D^2V(i\delta,\sqrt{\delta}S_i))\right)\right]\\
&=&L_1+L_2.
\end{eqnarray*}
Note that $\partial_{x_1x_2}^2V=\partial_{x_2x_1}^2V,$ so $D^2V\in\mathbb{S}(d)$, and thus $G(D^2V)$ is meaningful. Now we discuss $L_1$ and $L_2$ respectively. For $L_1$, by the definition of the function $G_{i+1}$, we have
\begin{eqnarray*}
L_1\le \frac{1}{n}\sum_{i=0}^{n-1}\hat{E}\left[\left(\frac{1}{2}\langle D^2V(i\delta,\sqrt{\delta}S_i)X_{i+1},X_{i+1}\rangle-
G_{i+1}(D^2V(i\delta,\sqrt{\delta}S_i))\right)\right]=0.
\end{eqnarray*}
For $L_2$, by condition (iv), we have
\begin{eqnarray*}
L_2\le \frac{1}{n}\hat{E}\left[\sum_{i=0}^{n-1}a_{i+1}\|D^2V(i\delta,
\sqrt{\delta}S_i)\|\right]    \le \frac{1}{n}\sum_{i=0}^{n-1}a_{i+1}\hat{E}\left[
\|D^2V(i\delta,\sqrt{\delta}S_i)\|\right].
\end{eqnarray*}
We claim that there exists a positive constant $C_1$ such that  $\hat{E}[\|D^2V(i\delta,\sqrt{\delta}S_i)\|]\le C_1,\forall i=1,2,\ldots$ By (\ref{lem4.1-c}), we have
$$
\|D^2V(i\delta ,\sqrt{\delta}S_i)-D^2V(0,0)\|\le C \left( |\sqrt{\delta}S_i|^{\alpha}+|i\delta|^{\frac{\alpha}{2}}\right),
$$
where $C$ is a positive constant. By H\"{o}lder's inequality and conditions (i)-(iii), we have
\begin{eqnarray*}
\hat{E}[|\sqrt{\delta}S_i|^{\alpha}]\le \left(\hat{E}[|\sqrt{\delta}S_i|^2]\right)^\frac{\alpha}{2}
=\left(\frac{1}{n}\sum _{k=1}^i\hat{E}[|X_k|^2]\right)^\frac{\alpha}{2}                                    \le\left(\frac{1}{n}\sum _{k=1}^i\left(\hat{E}[|X_k|^3]\right)^{\frac{2}{3}}\right)^\frac{\alpha}{2}                                    \le M^{\frac{2}{3}\times\frac{\alpha}{2}}                                    =M^{\frac{\alpha}{3}}.
\end{eqnarray*}
Let $C_1:=\|D^2V(0,0)\|+CM^{\frac{\alpha}{3}}+C$. Then for any $i=1,2,\ldots$,
$$
\hat{E}[\|D^2V(i\delta,\sqrt{\delta}S_i)\|]\le\|D^2V(0,0)\|
+C\hat{E}[|\sqrt{\delta}S_i|^{\alpha}]+C|i\delta|^{\frac{\alpha}{2}}
\leq C_1.
$$
Hence
$L_2\le \frac{C_1}{n}\sum_{i=1}^{n}a_{i}$, which together with condition (iv) implies that
$\limsup\limits_{n \to \infty}  L_2\le 0.$
Then
\begin{equation}\label{eq14}
\limsup\limits_{n \to \infty} \hat{E}\left[\sum _{i=0}^{n-1}J_{\delta}^{i}\right]\le 0.
\end{equation}
Similarly, we have
\begin{equation}\label{eq15}
\limsup\limits_{n \to \infty}\hat{E}\left[-\sum _{i=0}^{n-1}J_{\delta}^{i}\right]\le 0.
\end{equation}
As
$\hat{E}[\sum _{i=0}^{n-1}J_{\delta}^{i}]\ge -\hat{E}[-\sum _{i=0}^{n-1}J_{\delta}^{i}].$
By $(\ref{eq15})$, we have
$ \liminf\limits_{n \to \infty} \hat{E}[\sum _{i=0}^{n-1}J_{\delta}^{i}]\ge 0$, which together with  $(\ref{eq14})$ implies that
$
{\mathop {\lim }\limits_{n \to \infty } } \hat{E}\left[\sum _{i=0}^{n-1}J_{\delta}^{i}\right]= 0.
$

Now we come to analyze the two items including $I_{\delta}^i$ in (\ref{eq6}). Since $\partial _t V$ and $D^2V$ are uniformly $\alpha$-H\"{o}lder  continuous in $x$ and $\frac{\alpha}{2}$-H\"{o}lder continuous in $t$ on $[0,1]\times \mathbb{R}^d$, we have
$$|I_{\delta}^i|\le C^{'} \delta^{1+{\frac{\alpha}{2}}}(1+|X_{i+1}|^{\alpha}+|X_{i+1}
|^{2+\alpha}),$$
where $C'$ is a positive constant. By Young's inequality, we have
$$|X_{i+1}|^\alpha\le\frac{|X_{i+1}|^{\alpha\cdot\frac{2+\alpha}{\alpha}}}{\frac{2+\alpha}{\alpha}}+\frac{1}{\frac{2+\alpha}{2}}
                  =\frac{\alpha|X_{i+1}|^{2+\alpha}}{2+\alpha}+\frac{2}{2+\alpha},$$
and thus there exists a positive constant $C_2$ such that
$$|I_{\delta}^i|\le C_2 \delta ^{1+{\frac{\alpha}{2}}}(1+|X_{i+1}|^{2+\alpha}).$$
Hence
\begin{eqnarray*}
\hat{E}[|I_{\delta}^i|]&\le& C_2\delta^{1+\frac{\alpha}{2}}\left(1+\hat{E}[|X_{i+1}|^{2+\alpha}]\right)\\
&\le& C_2 \delta^{1+\frac{\alpha}{2}}\left(1+\left(\hat{E}
[|X_{i+1}|^3]\right)^{\frac{2+\alpha}{3}}\right)\\
&\le& C_2\delta^{1+\frac{\alpha}{2}}(1+M^{\frac{2+\alpha}{3}}).
\end{eqnarray*}
It follows that
$$\hat{E}\left[\sum _{i=0}^{n-1}I_{\delta}^i\right]\le \sum _{i=0}^{n-1}\hat{E}[I_{\delta}^i]\le \sum _{i=0}^{n-1}\hat{E}[|I_{\delta}^i|]\le
C_2(1+M^{\frac{2+\alpha}{3}})\left(\frac{1}{n}\right)^{\frac{\alpha}{2}},$$
and thus
\begin{equation}\label{eq13}
\limsup\limits_{n \to \infty} \hat{E}\left[\sum _{i=0}^{n-1}I_{\delta}^{i}\right]\le 0.
\end{equation}
On the other side,
$$\hat{E}\left[\sum _{i=0}^{n-1}I_{\delta}^i\right]\ge -\sum _{i=0}^{n-1}\hat{E}[-I_{\delta}^i]\ge -\sum _{i=0}^{n-1}\hat{E}[|I_{\delta}^i|]\ge
-C_2(1+M^{\frac{2+\alpha}{3}})\left(\frac{1}{n}\right)^{\frac{\alpha}{2}},$$
we have
$\liminf\limits_{n \to \infty} \hat{E}\left[\sum _{i=0}^{n-1}I^i_{\delta}\right]\ge 0,$
which together with $(\ref{eq13})$ implies that
$\lim\limits_{n \to \infty } \hat{E}[\sum _{i=0}^{n-1}I^i_{\delta}]=0.$
Similarly, we have
$
\lim \limits_{n \to \infty } \hat{E}[-\sum _{i=0}^{n-1}I^i_{\delta}]=0.
$

By $(\ref{eq6})$, we have
\begin{equation}\label{eq2}
\mathop {\lim }\limits_{n \to \infty } \hat{E}[V(1,\sqrt{\delta}S_n)]=V(0,0).
\end{equation}

Finally, $\forall t$, $t'\in [0,1+h]$, $x\in \mathbb{R}^d,$
\begin{eqnarray}\label{lem4.1-d}
|V(t,x)-V(t',x)|&\le&
\hat{E}[|\varphi(x+\sqrt{1+h-t}X)-\varphi(x+\sqrt{1+h-t'}X)|]\nonumber\\
                &\le& k_{\varphi}|\sqrt{1+h-t}-\sqrt{1+h-t'}|\hat{E}[|X|]\nonumber\\
                &\le& C_3 \sqrt{|t-t'|},
\end{eqnarray}
where $k_{\varphi}$ is a Lipschitz constant depending on $\varphi$, and $C_3$ is a constant depending on $\varphi$. By (\ref{lem4.1-b-1}) and (\ref{lem4.1-d}), we have
$$|\hat{E}[\varphi(\sqrt{1+h}X)]-\hat{E}[\varphi(X)]|=|V(h,0)-V(0,0)|\le C_3 \sqrt{h},$$
and
$$|\hat{E}\left[\varphi\left(\frac{S_n}{\sqrt{n}}\right)\right]-\hat{E}[V(1,\sqrt{\delta}S_n)]
|=|\hat{E}[V(1+h,\sqrt{\delta}S_n)]-\hat{E}[V(1,\sqrt{\delta}S_n)]|\le C_3\sqrt{h}.$$
It follows that
\begin{eqnarray*}
   \left|\hat{E}\left[\varphi\left(\frac{S_n}{\sqrt{n}}\right)\right]-
   \hat{E}[\varphi(X)]\right|
   &\le&\left|\hat{E}\left[\varphi\left(\frac{S_n}{\sqrt{n}}\right)\right]-
   \hat{E}[V(1,\sqrt{\delta}S_n)]\right|
   +|\hat{E}[V(1,\sqrt{\delta}S_n)]-\hat{E}[\varphi(\sqrt{1+h}X)]|\\
   &&+|\hat{E}[\varphi(\sqrt{1+h}X)]-\hat{E}[\varphi(X)]|\\
&\le& 2C_3\sqrt{h}+|\hat{E}[V(1,\sqrt{\delta}S_n)]-V(0,0)|,
\end{eqnarray*}
which together with (\ref{eq2}) implies that
$$
\limsup\limits_{n \to \infty } \left|\hat{E}\left[\varphi\left(\frac{S_n}{\sqrt{n}}\right)\right]
-\hat{E}[\varphi(X)]\right|\le 2C_3\sqrt{h}.
$$
Since $h$ can be arbitrarily small, we have
$$\mathop {\lim }\limits_{n \to \infty } \hat{E}\left[\varphi\left(\frac{S_n}{\sqrt{n}}\right)\right]=
\hat{E}[\varphi(X)].$$
\hfill\fbox

\bigskip

\noindent{\bf {Proof of Theorem \ref{Thm3.1}.}} The idea comes from the proof of Theorem 3.5 of \cite{Peng6}.
 When the uniform elliptic condition (\ref{lem4.1-a})
does not hold, we first introduce a perturbation to prove the above convergence
for $\varphi\in C_{b,Lip(\mathbb{R}^d)}$. We can construct a sublinear expectation space $(\bar{\Omega},\bar{\mathcal{H}},\bar{\mathbb{E}})$ and a sequence of two
random vectors $\{(\bar{X_i},\bar{\kappa_i})\}_{i=1}^{\infty}$ such that, for each $n=1,2,\ldots,\,\{\bar{X}_i\}_{i=1}^{n}\thicksim \{X_i\}_{i=1}^n$ and $\{(\bar{X}_{n+1},\bar{\kappa}_{n+1})\}$ is independent from $\{(\bar{X}_i,\bar{\kappa}_i)\}_{i=1}^{n}$ and,
moreover,
$$
\bar{\mathbb{E}}[\psi(\bar{X}_i,\bar{\kappa}_i)]
=(2\pi)^{-\frac{d}{2}}\int_{\mathbb{R}^d}\hat{E}[\psi(X_i,x)]
e^{-\frac{|x|^2}{2}}dx,~~\forall \psi \in C_{b,Lip}(\mathbb{R}^{2\times d}).
$$
Define $\bar{X}^{\varepsilon}_i:=\bar{X}_i+\varepsilon \bar{\kappa}_i$ for a fixed $\varepsilon >0$. It's easy to check that the sequence $\{\bar{X}^{\varepsilon}_i\}_{i=1}^{\infty}$ satisfies all the conditions of Theorem \ref{Thm3.1}, in particular,
$$G_i^{\varepsilon}(A):=\frac{1}{2}\bar{\mathbb{E}}[\langle A\bar{X}^{\varepsilon}_i,\bar{X}^{\varepsilon}_i\rangle]=G_i(A)
+\frac{\varepsilon ^2}{2}{\rm tr}[A],
$$
$$
G^{\varepsilon}(A):=G(A)+\frac{\varepsilon ^2}{2}{\rm tr}[A].
$$
They are strictly elliptic. Then we can apply  Lemma \ref{lem4.1} to
$$
\bar{S}_n^{\varepsilon}=\sum_{i=1}^n\bar{X}_i^{\varepsilon}=
\sum_{i=1}^n \bar{X}_i+\varepsilon J_n,~J_n=\sum_{i=1}^n \bar{\kappa}_i,
$$
and obtain
$$
\mathop {\lim }\limits_{n \to \infty } \bar{\mathbb{E}}\left[\varphi\left(\frac{\bar{S}_n^{\varepsilon}}
{\sqrt{n}}\right)\right]=\bar{\mathbb{E}}[\varphi(\bar{X}+\varepsilon \bar{\kappa})],
$$
where $(\bar{X},\bar{\kappa})$ is $\bar{G}$-distributed under $\bar{\mathbb{E}}[\cdot]$ and
$$
\bar{G}(\bar{A}):=\frac{1}{2}\bar{\mathbb{E}}[\langle \bar{A}(\bar{X},\bar{\kappa}),(\bar{X},\bar{\kappa})
\rangle],~~\bar{A}\in \mathbb{S}(2d).
$$
We have
\begin{eqnarray*}
\left|\hat{E}\left[\varphi\left(\frac{{S}_n}{\sqrt{n}}\right)\right]-
\bar{\mathbb{E}}\left[\varphi\left(\frac{\bar{S}_n^{\varepsilon}}
{\sqrt{n}}\right)\right]\right|
&=&\left|\bar{\mathbb{E}}\left[\varphi\left(\frac{\bar{S}_n^{\varepsilon}}
{\sqrt{n}}-\varepsilon \frac{J_n}{\sqrt{n}}\right)\right]-\bar{\mathbb{E}}\left[
\varphi\left(\frac{\bar{S}_n^{\varepsilon}}{\sqrt{n}}\right)\right]
\right|\\
&\le& \varepsilon C\bar{\mathbb{E}}\left[\left|\frac{J_n}{\sqrt{n}}\right|\right]\le C'\varepsilon,
\end{eqnarray*}
where $C$ is a Lipschitz constant of $\varphi$ and $C'$ is a constant depending on $\varphi$. Similarly,
$$
|\hat{E}[\varphi(X)]-\bar{\mathbb{E}}[\varphi(\bar{X}+\varepsilon \bar{\kappa})]|=
|\bar{\mathbb{E}}[\varphi(\bar{X})]-\bar{\mathbb{E}}[\varphi
(\bar{X}+\varepsilon \bar{\kappa})]|\le C \varepsilon.
$$
Since $\varepsilon$ can be arbitrarily small, it follows that
$$\mathop {\lim }\limits_{n \to \infty } \hat{E}\left[\varphi\left(\frac{S_n}{\sqrt{n}}\right)\right]=\hat{E}[\varphi(X)],~~\forall \varphi \in C_{b,Lip}(\mathbb{R}^d).$$
The proof of Theorem $\ref{Thm3.1}$ is complete.\hfill\fbox

\bigskip

\begin{remark}
By the proof of Lemma $\ref{lem4.1}$, we know that   condition (iii) in Theorem $\ref{Thm3.1}$ can be weaken to: there exist two constants $M>0$ and  $\alpha\in (0,1)$ such that $\hat{E}[|X_i|^{2+\alpha}]\le M,\forall i=1,2,\ldots$
\end{remark}

\subsection{Proof of Corollary \ref{Cor3.2}}

At first, we prove two lemmas.

\begin{lemma}\label{lem4.3}
We make the same assumptions as in Theorem \ref{Thm3.1}, then $\{G_i,G\}$ are Lipschitz function, and further they have a common Lipschitz constant $C_0$.
\end{lemma}
{\bf Proof.}
$\forall A,B \in \mathbb{S}(d)$ and $\forall i=1,2,\ldots$, we have
\begin{eqnarray*}
|G_i(A)-G_i(B)|& =&\frac{1}{2}|\hat{E}[\langle AX_i,X_i\rangle]-\hat{E}[\langle BX_i,X_i\rangle]|\\
                                 &\le& \frac{1}{2}\hat{E}[|\langle (A-B)X_i,X_i\rangle|]\\
                                 &\le&\frac{1}{2}\|A-B\|\cdot\hat{E}[|X_i|^2]\\
                                 &\le&\frac{1}{2}\|A-B\|\cdot\left(\hat{E}[|X_i|^3]\right)^{\frac{2}{3}}\\
                                 &\le&\frac{1}{2}M^{\frac{2}{3}}\|A-B\|,
\end{eqnarray*}
and
\begin{eqnarray*}
|G(A)-G(B)|&= &\lim_{n\rightarrow \infty}|G_i(A)-G_i(B)|
                             \le\frac{1}{2}M^{\frac{2}{3}}\|A-B\|.
\end{eqnarray*}
Set $C_0=\frac{1}{2}M^{\frac{2}{3}}$. The proof is complete.
\hfill\fbox

\begin{lemma}\label{lem4.4}
We make the same assumptions as in Theorem \ref{Thm3.1}, then there exists a  sequence $\left\{\gamma _i\right\}$ of positive numbers such that ${\mathop {\lim }\limits_{i \to \infty } }\gamma_i =0$ and
$$|G_i(A)-G(A)|\le \gamma _i \|A\|,\ \forall A\in \mathbb{S}(d).$$
\end{lemma}
{\bf Proof.}
Because of $G_i$'s and $G$'s positive homogeneity, we  only need to prove $\left\{G_i \right\}$ uniformly converges to $G$ on the set $\{A\in \mathbb{S}(d)|\,\|A\|\le 1\}$. We suppose that this does not hold. Then there exists an $\varepsilon >0$ such that $\forall k$, $\exists~ n_k\ge k$, $A_{n_k}\in\{A\in \mathbb{S}(d)\mid \|A\|\le 1\}$, and
$$
|G_{n_k}(A_{n_k})-G(A_{n_k})|>\varepsilon.
$$
 Since $\{A\in \mathbb{S}(d)\mid ||A||\le 1\}$  is a compact subset of $\mathbb{S}(d)$, thus it is sequential compact, so $\{A_{n_k}\}$ has a convergent subsequence $\{A_{n_{k_l}}\}$. Denote by $A_0$ its limit point. By Lemma \ref{lem4.3}, we get
$$|G_{n_{k_l}}(A_{n_{k_l}})-G_{n_{k_l}}(A_0)|\le C_0\|A_{n_{k_l}}-A_0\|,\ |G(A_{n_{k_l}})-G(A_0)|\le C_0\|A_{n_{k_l}}-A_0\|,$$
and thus
$$|G_{n_{k_l}}(A_{n_{k_l}})-G(A_{n_{k_l}})|\le |G_{n_{k_l}}(A_{n_{k_l}})-G_{n_{k_l}}(A_0)|+|G_{n_{k_l}}(A_0)-G(A_0)|+
|G(A_0)-G(A_{n_{k_l}})|.$$
Leting $l\to \infty$,  we obtain that $\limsup\limits_{l \to \infty }  |G_{n_{k_l}}(A_{n_{k_l}})-G(A_{n_{k_l}})|\le 0$. It's a contradiction. Hence the result of this lemma holds.
\hfill\fbox

\smallskip

Then  Corollary \ref{Cor3.2} directly follows from Theorem \ref{Thm3.1} and Lemma \ref{lem4.4}.

\subsection{Proof of Corollary \ref{cor3.4}}

Set $$G(A)=\frac{1}{2}\hat{E}[\langle AX,X\rangle]=\frac{1}{2}\sup_{Q\in\Theta}{\rm tr}[AQ],\forall A\in \mathbb{S}(d).
$$
We claim that
\begin{eqnarray}\label{pf-cor3.4-a}
|G_i(A)-G(A)|\le d_H(\Theta_i,\Theta)\|A\|.
\end{eqnarray}
If $d_H(\Theta_i,\Theta)= 0$, then $\Theta_i=\Theta$, and so $G_i=G$ and thus (\ref{pf-cor3.4-a}) holds in this case. Now assume that $d_H(\Theta_i,\Theta)\neq 0$.
By ${\rm tr}[AQ_1]-{\rm tr}[AQ_2]={\rm tr}[A(Q_1-Q_2)],$ we have
\begin{equation}\label{pf-cor3.4-b}
|{\rm tr}[AQ_1]-{\rm tr}[AQ_2]|\le\|A\|\cdot\|Q_1-Q_2\|.
\end{equation}
Suppose that (\ref{pf-cor3.4-a}) doesn't hold. Then
there exists an element $A\ne 0\in \mathbb{S}(d)$ such that
\begin{eqnarray}\label{pf-cor3.4-c}
|G_i(A)-G(A)|> d_H(\Theta_i,\Theta)\|A\|.
\end{eqnarray}
 Since $\Theta_i$ and $\Theta$ are bounded and closed,  there exist $Q_i\in\Theta_i, Q\in\Theta$ such that
 \begin{eqnarray}\label{pf-cor3.4-c-1}
 G_i(A)=\frac{1}{2}{\rm tr}[AQ_i],\ \ G(A)=\frac{1}{2}{\rm tr}[AQ].
 \end{eqnarray}
Without loss of generality, we assume that $G_i(A)\ge G(A)$.  Then by (\ref{pf-cor3.4-c}) and (\ref{pf-cor3.4-c-1}), we have
\begin{eqnarray}\label{pf-cor3.4-d}
\frac{1}{2}{\rm tr}[AQ_i]-\frac{1}{2}{\rm tr}[AQ]=\frac{1}{2}{\rm tr}[A(Q_i-Q)]> d_H(\Theta_i,\Theta)\|A\|.
\end{eqnarray}
By the definition of $d_H$ and the assumption that $d_H(\Theta_i,\Theta)\neq 0$, there exists $Q_0\in\Theta$ such that \begin{eqnarray}\label{pf-cor3.4-e}
\|Q_i-Q_0\|<\frac{3}{2}d_H(\Theta_i,\Theta).
\end{eqnarray}
Then by (\ref{pf-cor3.4-b}), (\ref{pf-cor3.4-d}) and (\ref{pf-cor3.4-e}), we have
\begin{eqnarray*}
\frac{1}{2}{\rm tr}[AQ_0]-\frac{1}{2}{\rm tr}[AQ]
&=&\frac{1}{2}{\rm tr}[AQ_0]-\frac{1}{2}{\rm tr}[AQ_i]+\frac{1}{2}{\rm tr}[A(Q_i-Q)]\\
&>&-\frac{1}{2}\|A\|\cdot\|Q_0-Q_i\|+d_H(\Theta_i,\Theta)\|A\|\\
&\ge&-\frac{3}{4}d_H(\Theta_i,\Theta)\|A\|+d_H(\Theta_i,\Theta)\|A\|\\
& =&\frac{1}{4}d_H(\Theta_i,\Theta)\|A\|>0,
\end{eqnarray*}
which contradicts the fact that $\frac{1}{2}{\rm tr}[AQ]=G(A)=\frac{1}{2}\sup_{Q\in\Theta}{\rm tr}[AQ]$. Hence (\ref{pf-cor3.4-a})  holds. Then by Theorem \ref{Thm3.1}, we obtain the result of Corollary \ref{cor3.4}.\hfill\fbox

\subsection{Proof of Corollary \ref{Cor3.5}}

Since $\{\Theta _i\}$ is a Cauchy sequence in $(\Xi,d_H)$ and $(\Xi,d_H)$ is a complete metric space, there exists a bounded closed set $\Theta\in\Xi$ such that
$\lim\limits_{n\to \infty}d_H(\Theta,\Theta_i)=0,$ hence
$\lim\limits_{n\to \infty}\frac{1}{n}\sum_{i=1}^{n}d_H(\Theta,\Theta_i)=0.$
And thus by Corollary \ref{cor3.4}, we only need to  prove that $\Theta$ is a convex set. First, we prove two claims.

{\it Claim 1.} $x\in\Theta\Leftrightarrow\lim\limits_{i\to \infty}d(x,\Theta_i)=0.$

``$\Rightarrow$" If $x\in\Theta$, then by the definition of $d_H$, we have $d(x,\Theta_i)\le d_H(\Theta,\Theta_i),$ and
thus $\lim\limits_{i\rightarrow \infty}d(x,\Theta_i)=0.$

``$\Leftarrow$" Since $d_H(\Theta,\Theta_i)\to 0$ as $i\rightarrow \infty,$ so $\forall \varepsilon>0$, there exists an $N$ such that when $i>N$, $d_H(\Theta,\Theta_i)<\varepsilon/2.$
Hence $\Theta_i\subset B_\epsilon(\Theta)$, and by $\lim\limits_{i\rightarrow \infty}d(x,\Theta_i)=0,$
 we get $x\in B_{2\epsilon}(\Theta)$. Since $\varepsilon$ is can be arbitrarily small, we have $x\in\overline{\Theta}=\Theta.$

{\it Claim 2.} The function $d(\cdot,\Theta_i):\mathbb{S}_+(d)\to \mathbb{R}$ is convex.

Indeed, $\forall x,y\in\mathbb{S}_+(d)$, $a\in[0,1]$, there exist $x_0,\;y_0\in\Theta_i$ such that
$$d(x,\Theta_i)=d(x,x_0),~d(y,\Theta_i)=d(y,y_0).$$
Since $\Theta_i$ is convex, we have $ax_0+(1-a)y_0\in\Theta_i$.
Then we have
\begin{eqnarray*}
d(ax+(1-a)y,\Theta_i)
&\le&d(ax+(1-a)y,ax_0+(1-a)y_0)\\
& =&\|\left(ax+(1-a)y\right)-\left(ax_0+(1-a)y_0\right)\|\\
&\le&a\|x-x_0\|+(1-a)\|y-y_0\|\\
& =&ad(x,\Theta_i)+(1-a)d(y,\Theta_i),
\end{eqnarray*}
and thus  {\it Claim 2} holds.

Now we prove  that $\Theta$ is a convex set.
For any  $x,y\in\Theta$ and $a\in[0,1]$, by {\it Claim 1} and {\it Claim 2},  we have
$$\lim_{i\to\infty}d(ax+(1-a)y,\Theta_i)\le\lim_{i\to\infty}
(ad(x,\Theta_i)+(1-a)d(y,\Theta_i))=0.$$
By {\it Claim 1} again, we get that $ax+(1-a)y\in\Theta$, and
thus $\Theta$ is a convex set. \hfill\fbox

\bigskip

\begin{remark}
Under the conditions of Corollary \ref{Cor3.5}, by the above  proof and   (\ref{pf-cor3.4-a}) we have $\lim\limits_{i\to \infty}G_i(A)=G(A), \forall A\in \mathbb{S}(d).$
\end{remark}

\smallskip

\section{Multi-dimensional laws of large numbers under sublinear expectations}

In this section, we give a  multi-dimensional (weak) law of large numbers under sublinear expectations and several corollaries. The idea for the proofs is similar to
the one used in the above section and  we omit the details.

\begin{definition}(\cite{Peng6})
A $d$-dimensional random vector $\eta=(\eta_1,\eta_2,\ldots,\eta_d)$ on a sublinear expectation space $(\Omega,\mathcal{H},\hat{E})$ is called maximal distributed if there exists a bounded, vonvex and closed subset $\Gamma\in\mathbb{R}^d$ such that
$$\hat{E}[\varphi(\eta)]=\max_{y\in\Gamma}\varphi(y),~\forall \varphi\in C_{b,Lip}(\mathbb{R}^d).$$
\end{definition}

\smallskip

\begin{proposition}(\cite{Peng6})\label{Pro4}
Let $\eta$ be maximal distributed. For any $\varphi \in C_{b,Lip}(\mathbb{R}^d)$,
define a function
$$u(t,y)=\hat{E}[\varphi(y+t\eta)],~(t,y)\in [0,\infty)\times \mathbb{R}^d.$$
Then u is the unique viscosity solution of the following parabolic PDE:
 \begin{equation}\label{eq23}
\partial_t u-g(Du)=0,~~u|_{t=0}=\varphi,
\end{equation}
where
$g=g_{\eta}:\mathbb{R}^d\to\mathbb{R}$ is defined by $g_{\eta}(p)=\hat{E}[\langle p,\eta\rangle].$
\end{proposition}

\smallskip

\begin{remark}
It is easy to check that $g_{\eta}$ is a sublinear function on $\mathbb{R}^d$. Then by Theorem \ref{Th1}, there exists a bounded, convex and closed subset $\Gamma\subset \mathbb{R}^d$ such that
$$g(p)=\sup_{q\in \Gamma}\langle p,q\rangle,~~\forall p\in \mathbb{R}^d.$$
\end{remark}

\smallskip

\begin{theorem}\label{Thm5.4}
Let  $\{Y_i\}_{i=1}^{\infty}$ be a sequence of $\ \mathbb{R}^d$-valued random variables
in a sublinear expectation space $(\Omega,\mathcal{H},\hat{E})$ which satisfies the following conditions:
\begin{itemize}
\item[(i)] each $Y_{i+1}$ is independent from $(Y_1,\ldots,Y_i)$, $\forall i=1,2,\ldots;$
\item[(ii)] there exists a positive constant $M$ such that $\hat{E}[|Y_i|^2]\le M,\,\forall i=1,2,\ldots;$
 \item[(iii)] there exist a  sequence $\{a_n\}$ of positive numbers and a sub-additive, positive homogeneous function ${g}:\mathbb{R}^d\rightarrow\mathbb{R}$ such that
$\lim\limits_{n\to \infty}(a_1+\cdots+a_n)/n=0,$ and
\begin{eqnarray*}
|{g}_n(p)- {g}(p)|\le a_n|p|,\forall n,~~p\in \mathbb{R}^d,
\end{eqnarray*}
where ${g}_i:\mathbb{R}^d\rightarrow\mathbb{R},~{g}_i(p)=\hat{E}[\langle Y_i,p\rangle],~~\forall p\in\mathbb{R}^{d}.$
\end{itemize}
Then the sequence $\left\{\frac{S_n}{{n}}\right\}_{n=1}^{\infty}$ converges in law to G-normal distribution $N(\Gamma;0)$, i.e.,
$$\lim_{n\rightarrow\infty}\hat{E}\left[\varphi\left(\frac{S_n}{{n}}\right)\right]=\hat{E}[\varphi(\eta)],~\varphi\in C_{b,Lip}(\mathbb{R}^d),$$
where  $S_n=Y_1+\cdots+Y_n$, $\eta\thicksim N(\Gamma;0)$ and $g(p)=\sup_{q\in \Gamma}\langle p,q\rangle,~~\forall p\in \mathbb{R}^d.$
\end{theorem}

\smallskip

\begin{corollary}
Let  $\{Y_i\}_{i=1}^{\infty}$ be a sequence of $\ \mathbb{R}^d$-valued random variables
in a sublinear expectation space $(\Omega,\mathcal{H},\hat{E})$ which satisfies  conditions (i)(ii) in Theorem \ref{Thm5.4} and the following condition (iv).
\begin{itemize}
\item[(iv)] $\{g_i\}$  converges to ${g}$ pointwise, i.e.,
$\mathop {\lim }\limits_{i \to \infty } g_i(p)=g(p),~\forall p\in \mathbb{R}^d,$ where ${g}_i:\mathbb{R}^d\rightarrow\mathbb{R},~{g}_i(p)=\hat{E}[\langle Y_i,p\rangle],~\forall p\in\mathbb{R}.$
\end{itemize}
Then the result of Theorem \ref{Thm5.4} holds.
\end{corollary}

\smallskip

Define
$$\Pi:=\{\Gamma\subset\mathbb{R}^d| \Gamma \ \mbox{is  bounded and closed}\}.$$
 Denote by $d_H$ the Hausdorff distance on $\Pi$. Then
by \cite{He}, $(\Pi,d_H)$ is a complete metric space.

\smallskip

\begin{corollary}\label{cor5.6}
Let  $\{Y_i\}_{i=1}^{\infty}$ be a sequence of $\ \mathbb{R}^d$-valued random variables
in a sublinear expectation space $(\Omega,\mathcal{H},\hat{E})$ which satisfies  conditions (i)(ii) in Theorem \ref{Thm5.4} and the following condition (v).
\begin{itemize}
\item[(v)] There exists a bounded, convex and  closed  subset $\Gamma\subset\mathbb{R}^d$ such that
    $$\lim _{n\to \infty}\frac{1}{n}\sum_{i=1}^{n}d_H(\Gamma_i,\Gamma)=0,$$
where $\Gamma_i$ is a bounded, convex and closed  subset of $\ \mathbb{R}^d$ such that $\hat{E}[\langle Y_i,p\rangle]=
\sup_{q\in \Gamma_i\subset \mathbb{R}^d}\langle q,p\rangle$, $\forall p\in\mathbb{R}^d.$
\end{itemize}
Then the result of Theorem \ref{Thm5.4} holds.
\end{corollary}

\smallskip

\begin{corollary}
Let  $\{Y_i\}_{i=1}^{\infty}$ be a sequence of $\ \mathbb{R}^d$-valued random variables
in a sublinear expectation space $(\Omega,\mathcal{H},\hat{E})$ which satisfies  conditions (i)(ii) in Theorem \ref{Thm5.4} and the following condition (vi).
\begin{itemize}
\item[(vi)] Let $\Gamma_i$ be a bounded, convex and  closed subset of $\ \mathbb{R}^d$ such that $\hat{E}[\langle Y_i,p\rangle]=
\sup_{q\in \Gamma_i\subset \mathbb{R}^d}\langle q,p\rangle,\forall p\in\mathbb{R}^d.$ Suppose that $\{\Gamma_i\}$ is a Cauchy sequence in $(\Pi,d_H)$.
\end{itemize}
Then there exists a bounded, convex and  closed  subset $\Gamma$ of $\ \mathbb{R}^d$ such that $d_H(\Gamma,\Gamma_i)\to 0$ as $i\to \infty$ and  the result of Theorem \ref{Thm5.4} holds.
\end{corollary}

\begin{remark}
By Corollary \ref{cor5.6} and the definition of the Hausdorff metric $d_H$, following the proof of Theorem \ref{thm3.6}, we can obtain Theorem 3.1 of \cite{HuF}.
\end{remark}

%

\end{document}